\documentclass[leqno,11pt]{article}
\usepackage{amssymb,amsmath,amscd, mathrsfs}

\makeatletter
\@addtoreset{equation}{subsection} 

\renewcommand{\@begintheorem}[2]{
\begin{trivlist}\it \item[\hspace{\labelsep}{\bf #1\ #2\ }]}

\renewcommand{\@opargbegintheorem}[3]{
\begin{trivlist}\it \item[\hspace{\labelsep}{\bf #1\ #2\ (#3)\ }]}

\renewcommand{\@endtheorem}{\end{trivlist}}
\makeatother

\newtheorem{teo}[equation]{Theorem}

\newtheorem{lem}[equation]{Lemma}

\newtheorem{pro}[equation]{Proposition}
\newtheorem{cor}[equation]{Corollary}

\newtheorem{teo1}{Theorem}

\newenvironment{ea*}{\begin{eqnarray*}}{\end{eqnarray*}}

\newenvironment{rem}{ \refstepcounter{equation}\par\noindent
 {{\bf Remark\/}  \bf \thesubsection.\arabic{equation}\ \ }}{\par}

\def\AA{\mathcal{A}}

\def\LL{\mathcal{L}}

\def\CT{\widetilde{C}}
\def\NT{\widetilde{N}}

\def\THT{\widetilde{\Theta}}
\def\PR{\PP(\tilde{C},C,\Xi)}

\def\OO{\mathcal{O}}

\def\lra{\longrightarrow}

\def\PRRRR{\mathbb{P}^{4}}

\newcommand{\be}{\begin{enumerate}}
\newcommand{\MF}[1]{\mathfrak{#1}}

\renewcommand{\PR}[1]{\mathbb{P}^{#1}}

\def\si{\sigma}
\def\tt{\tau}

\def\CX{\mathbb{C}}
\def\ZZ{\mathbb{Z}}

\def\PIC{\textnormal{Pic}}

\def\SPEC{\textnormal{Spec }}
\def\MUL{\textnormal{mult}}

\def\QED{\hspace*{\stretch{1}$\square$}}
\def\PRF{\noindent\emph{Proof}}
\def\remskip{\vskip .12in}

\def\sing{\textnormal{Sing }}
\def\singd{\textnormal{Sing}}
\def\supp{\textnormal{supp}}

\def\NM{\textnormal{Nm}}
\newcommand{\comment}[1]{}

\begin{document}
\bibliographystyle{amsplain}

\title{Cubic threefolds and abelian varieties of dimension five. II}
\author{Sebastian Casalaina-Martin\thanks{The author was partially
supported by NSF MSPRF grant DMS-0503228}}

\maketitle

\begin{abstract}   This paper extends joint work with R. Friedman to show that the closure of
the locus of intermediate Jacobians of smooth cubic threefolds, in the moduli
space of principally polarized abelian varieties (ppavs) of dimension five, is
an irreducible component of the locus of ppavs whose theta divisor has a point
of multiplicity three or more. This paper also gives a sharp bound on the
multiplicity of a point on the theta divisor of an indecomposable ppav of
dimension less than or equal to five; for dimensions four and five, this
improves the bound due to J. Koll\'ar, R. Smith-R. Varley, and L. Ein-R.
Lazarsfeld. \end{abstract}

\section*{Introduction}

The geometric Schottky problem
 is to identify Jacobian varieties among all principally polarized abelian varieties
(ppavs) via geometric conditions on the polarization.   For a smooth cubic hypersurface  $X\subset \mathbb{P}^4$, which we will simply call a cubic threefold, the intermediate Jacobian $(JX,\Theta_X)$ is a ppav of dimension five  
and one can consider the analogous problem, the ``geometric
Schottky problem for cubic threefolds,''  which is to identify these intermediate Jacobians among all ppavs of dimension five via
geometric conditions on the polarization.

Since a theorem of Mumford's \cite{mum} states that $\Theta_X$ has a unique singularity, which is of multiplicity three, it is natural to ask to what extent the existence of triple points on the theta divisor of a ppav of dimension five characterizes intermediate Jacobians of cubic threefolds.  In fact, in \cite{cmf}, Friedman and the author showed that a
ppav of dimension five whose theta divisor has a unique singular point of
multiplicity three, and no other singular points,
 is the intermediate Jacobian of a
 cubic threefold.
 
The first goal of this paper is to extend this result:

\begin{teo1}\label{teo0} Let $I\subset \AA_5$ be the locus of 
intermediate Jacobians of  cubic threefolds, let $N_1$ be the
locus of ppavs in $\AA_5$ whose theta divisor has a singular locus of dimension at
least one, and let $S_3$ be the locus of ppavs in $\AA_5$ whose theta divisor has a
point of multiplicity three or more. 
$$ I=S_3 - (N_1 \cap S_3).$$
\end{teo1}

\comment{
\begin{teo1}

The closure in $\AA_5$ of the locus of intermediate Jacobians of smooth cubic
threefolds in $\PR{4}$ is an irreducible component of the locus of ppavs whose
theta divisor has a point of multiplicity three or more. \end{teo1}
}
 
A more detailed statement is given in Theorem \ref{teos3}, which provides a complete
description of the locus $S_3$.  Roughly speaking, $S_3$ consists of three irreducible components, all of dimension ten, one of which is the closure of $I$.  
Theorem \ref{teos3} also gives new information on the boundary of $I$ (cf. Remark \ref{rembound}).


\remskip 
The techniques of this paper can also be used to give a sharp bound on the
multiplicity of a point on the theta divisor of an indecomposable ppav of dimension at
most five; for dimensions four and five, this improves the bound due to Koll\'ar
\cite{kollar},
 Smith-Varley \cite{sv4}, and Ein-Lazarsfeld \cite{einlaz}.  To be precise, let
$\textnormal{Sing}_k\Theta=\{x\in \Theta : \textnormal{mult}_x\Theta\ge k\}$.  A
result of Koll\'ar's \cite{kollar} shows that if $(A,\Theta)\in \mathcal{A}_d$, then
$\dim(\textnormal{Sing}_k\Theta)\le d-k$; generalizing a result of Smith and Varley
\cite{sv4}, Ein and Lazarsfeld \cite{einlaz} showed that
$\dim(\textnormal{Sing}_k\Theta) =d-k$ only if $(A,\Theta)$ splits as a $k$-fold
product. It follows that if $(A,\Theta)$ is indecomposable then $$
\dim(\textnormal{Sing}_k\Theta) \le d-k-1, $$ As a special case we see that if $x\in
\Theta$, then $\MUL_x\Theta\le d-1$. For $d\le 3$ it is easy to see that these
bounds are sharp; this paper shows that these bounds are not sharp for $d=4,5$. 

Motivation for this result comes from the case of Jacobian and Prym varieties.  For the
Jacobian $(JC,\Theta_C)$ of a smooth curve $C$ of genus $g$, applying the Riemann
singularity theorem and Martens' theorem \cite{martens}, it follows that
$\dim(\textnormal{Sing}_k\Theta_C)\le g-2k+1$, with equality holding only if $C$ is
hyperelliptic.  Similarly, in the case of an indecomposable Prym variety $(P,\Xi)$
associated to a connected \'etale double cover of a smooth curve $C$ of genus $g$,
the results of \cite{casa} show that $\dim(\textnormal{Sing}_k\Xi)\le (g-1)-2k+1$. 

The following theorem extends these results to all ppavs of dimension less than or
equal to five. 

\begin{teo1}\label{teoel} Suppose $(A,\Theta)\in \mathcal{A}_d$ for $d\le 5$. For
$k$ and $j$ nonnegative integers such that $k\ge 1$ and $k-j\ge1$, if $\dim(\textnormal{Sing}_k\Theta)= d-k-j$, then 
$(A,\Theta)=\prod_{i=1}^{k-j}(A_i,\Theta_i)$. Equivalently, if $(A,\Theta)$ is
indecomposable, $\dim(\textnormal{Sing}_k\Theta)\le d-2k+1$. \end{teo1}

This theorem is obtained from the following: 

\begin{teo1}\label{theorem}
 Suppose $(A,\Theta)\in \AA_d$ is indecomposable
 and either $d\le 5$ or  $(A,\Theta)$ is the Prym variety associated to the double cover of an irreducible stable curve.  For any $x\in \Theta$ we have $\MUL_x\Theta \le (d+1)/2$.
\end{teo1}

In this paper we will work over $\mathbb C$.  The techniques used here are similar to those used in \cite{cmf} and
\cite{casa}:  we use the fact due to Beauville \cite{bshot} that a ppav of dimension
at most five is the Prym variety associated to an admissible double cover of a
stable curve, and then, by studying deformations of line bundles on curves,
 describe the possible singularities of the theta divisor of a Prym variety. 
Section 1 recalls the basic setup in \cite{cmf} and \cite{casa}; a more general
situation is considered in this paper, and the technical difficulties this
introduces are dealt with in Section 1.2, culminating in Propositions \ref{mainpro1}
and \ref{mainpro2}. In Section 2 we apply the results of
Section 1 to the Prym theta divisor, and to description of $S_3\subseteq \mathcal{A}_5$.  The main theorems are proven in Sections 2.2 and 2.3. 

\subsection*{Acknowledgments}
It is a pleasure to thank Robert Friedman, as well as
Roy Smith and Robert Varley for many
insightful comments on Prym varieties and cubic threefolds.  I would also like to thank Samuel Grushevsky for his discussions on singularities of theta divisors, and the referee for many improvements in exposition.

\section{Preliminaries}\label{secpre}

\subsection{Prym varieties of nodal curves}

In this section, following Beauville \cite{bshot},
 we recall some basic
results about the theta divisors of Prym varieties. Throughout the paper, when discussing Pryms,
$\CT$ will be a connected curve with at worst ordinary double points, equipped with an involution $\tt:\CT \to \CT$ satisfying Beauville's condition
($\ast$): 

\begin{itemize}
\item[($\ast$)]
the fixed points of $\tt$ are exactly the singular points of $\CT$,
and at a singular point, the two branches are not exchanged under $\tt$.
\end{itemize}
It is easy to show that $C=\CT/(\tt)$ has only ordinary double points 
(\cite{bshot}, Lemma 3.1), and we will define $\pi:\CT \to C$  
to be the induced morphism.  Observe as well that under these conditions
on $\CT$, $\deg(\omega_{\CT} |_{\CT_i})$ is even for each irreducible
component $\CT_i\subseteq \CT$, so that there is a well-defined theta divisor
for $J\CT$.

\comment{
Slightly abusing notation, we will let  
$\nu:\NT \to \CT$, and $\nu: N\to C$ be the normalizations,
so that we
have a commutative diagram
\[
\begin{CD}
\NT   @>\nu>>   \CT\\ 
@VV\pi V                  @VV\pi V\\ 
N          @>\nu>>        C
\end{CD}
\]

Given a double cover of curves, we will let $\NM$ be the usual norm for
line bundles ([EGA II.6.5]). 
In the situation we have been considering, this induces a
commutative diagram

\begin{equation}\label{cdnorm}
\begin{CD}
0 @>>> \widetilde{T} @>>> J\CT   @>\nu^*>>   J\NT @>>> 0\\ 
@.  @VV \NM V                  @VV \NM V  @VV \NM V \\ 
0 @>>> T @>>> JC          @>\nu^*>>   JN @>>> 0,
\end{CD}
\end{equation}
where both $\widetilde{T}$ and $T$ are isomorphic to 
$(\CX^*)^{\delta-\mu+1}$.
We define the Prym variety $P\subseteq \ker (\NM)\subseteq J\CT$ 
to be the connected 
component of the identity; Beauville has shown that 
$\ker(\NM)\cong P \times \ZZ/2\ZZ$.  Also, since 
$\pi^*$ induces an isomorphism of $T$ onto $\widetilde{T}$, and 
since $\NM \circ \pi^*=2$, it follows that $\NM|_{\widetilde{T}}$ is surjective,
and $\ker(\NM|_{\widetilde{T}}) = \widetilde{T}_2 
\cong (\ZZ/2 \ZZ)^{\delta-\mu+1}$, the points of order $2$ in $\widetilde{T}$.
Thus we have an exact sequence,

$$
0 \lra \widetilde{T}_2 \lra P \times \ZZ/2\ZZ \lra R \lra 0,
$$
where $R=\ker(\NM|_{J\NT})$.  In particular, $P$ is isogenous to $R$.}

Let $\nu:\NT \to \CT$, and $\nu: N\to C$ be the normalization of the respective curves, and let
$\NM:\PIC(\CT)\to \PIC(C)$ be the usual norm for line bundles ([EGA II.6.5]).  In
analogy with the smooth case, the Prym variety $P\subseteq \ker (\NM)\subseteq J\CT$
is defined to be the connected component of the identity. Setting $JC^*$  to be the variety of line bundles $L$ on
$C$ such that $2c_1(L)=c_1(\omega_C)$, Beauville has shown that
$P$ can be identified with the set $$P^*=\{L\in J\CT^*\ | \ \NM(L)=\omega_C, \
h^0(L)\equiv 0 \ (\textnormal{mod 2})\}.$$ A principal polarization is given on $P$
by a divisor $\Xi$, which can be identified with the set $\Xi^*=\{L\in P^* \ |\
h^0(L)\ge 2\}$, and $\THT|_P=2\cdot \Xi$, where $\THT$ is the canonical polarization
of $J\CT$.  Thus, if $x\in P$ corresponds to a line bundle $L\in P^*$, then
$\MUL_x\Xi \ge \frac{1}{2}h^0(L)$.

\remskip Our starting point for studying $\Xi$ will be the following.
 Let $\LL$ be a line bundle over $C
\times J\CT$ of multidegree $\frac{1}{2}c_1(\omega_{\CT})$ with the property that for all
$x\in J\CT$, $\LL|_{\CT\times x}$ is isomorphic to the line bundle $M\in J\CT^*$ corresponding to
$x$. The existence of such a line bundle is established in Beauville \cite{bshot},
Lemma 1.3.

Now suppose that $S$ is a smooth curve with $s_0\in S$.  Let $t$ be a local coordinate for $S$ centered at $s_0$ and set $S_k=\SPEC \mathbb C[t]/(t^{k+1})$.  For each $k$ there is a map $S_k \to S$, so that if we set $\CT_k=\CT\times S_k$,  there are induced maps $\CT_k\to \CT\times S$.  Let $\mathcal L_k$ be the pull-back of $\LL$ to $\CT_k$.

The following was proven in \cite{cmf}.

\begin{lem}\label{lemest}
Suppose $S$ is a smooth curve, $f:S\to P$ is a morphism, and $f(s_0)=x\in P$
corresponds to a line bundle $L\in P^*$.
Let $\pi_2:C\times S\to S$ be the second projection, let $(\textnormal{id}_C\times f):C\times S\to C\times P$ be the induced morphism, and set 
$\LL'=(\textnormal{id}_C\times f)^*\LL$.  Then
$$
\frac{1}{2}h^0(L)\le \MUL_x\Xi \le \frac{1}{2}\deg_{s_0}\THT|_S
=\frac{1}{2}\ell((R^1\pi_{2*}\LL')_{s_0}),
$$
where $\ell((R^1\pi_{2*}\LL')_{s_0})$ is the length of $(R^1\pi_{2*}\LL')_{s_0}$ as an $\mathcal O_{S,s_0}$-module.
Moreover, there exist $S$, $f$, and $s_0$ as above such that $\MUL_x \Xi =
\frac{1}{2}\ell((R^1\pi_{2*}\LL')_{s_0})$. 
If $f(S)\nsubseteq \Xi$,  there is an $N\in \ZZ$ such that $\forall k\ge N$,
$$\ell((R^1\pi_{2*}\LL')_{s_0})=\ell(H^0(\LL_k)).$$
In particular, if for some $n\in \mathbb Z$, $\ell(H^0(\LL_n))=\ell(H^0(\LL_{n+1}))$, we may take $N=n$.
\QED
\end{lem}

\subsection{Line bundles on double covers}

The results of \cite{casa} Section 2 will be extended
to the case of nodal curves.  In particular, we will consider 
the case
where the double cover
$\pi:\CT\to C$ is a finite degree-two morphism of stable curves, induced by 
an involution $\tt:\CT \to \CT$, such that:

\begin{itemize}
\item[$(\S)$]the singular points are 
fixed points, and at the singular points, the two branches are not exchanged
under $\tt$.
\end{itemize}
This is  a slightly weaker condition than $(\ast)$, which
required that the fixed points be exactly the singular points.

We will begin by studying how points $p,\tt(p) \in \CT$ impose conditions on 
global sections of line bundles on $\CT$.  The difficulty here
is due to the fact that $\CT$ may be reducible, and hence nonzero 
sections may
vanish identically on
irreducible components.
We will fix the following notation: let $\nu:\NT \to \CT$ 
and $N\to C$ 
be the normalizations, let $\CT=\bigcup_{i=1}^\mu \CT_i$, 
where the $\CT_i$ are the irreducible components
of $\CT$, let $\sing\CT=\{x_1,\ldots,x_\delta \}$, let
$\nu^{-1}(x_i)=\{r_i,r_i'\}$, and let 
$\hat{R}=\sum_{i=1}^\delta (r_i+r_i')$.  
Recall that we are allowing
the covering to be ramified 
at some smooth points, say $\{t_1,\ldots,t_a\}$ of $\CT$; let
$\hat{T}=\sum_{i=1}^a t_i$. 

\remskip

The following definition generalizes the notation used in \cite{cmf}
and \cite{casa} to this setting:
for a line bundle $L$ on $\CT$,
a decomposition of $L$ is an 
isomorphism  
$$\phi: \nu^*L \to \pi^*M \otimes \OO_{\NT}(R +T+B),$$
where $M$ is a line bundle on $N$ such that
 $h^0(N,M)>0$,
 and $B$, $T$ and $R$ are effective divisors on $\NT$ satisfying the properties
 $B\cap \tau^*B=\varnothing$, 
 $T \le \hat{T}$,
$R \le \hat{R}$ and  $r_i \in \supp (R)$ if and only
if $r_i' \in \supp (R)$.

Given a decomposition of  $L$,
we will let
$r\in H^0(\NT,R)$, (resp. 
$t\in H^0(\NT, T)$, $b \in H^0(\NT, B)$), 
be a section vanishing on $R$, (resp. $T$, $B$).
We can consider $H^0(\CT,L)\subseteq H^0(\NT,\nu^*L)$ as a linear subspace;
a priori, the sections of $\pi^*H^0(N,M)\cdot r \cdot t\cdot b$ are not
necessarily sections of $H^0(\CT, L)$.  Define 
$\MF{M} \subseteq H^0(N,M)$ to be the maximal subspace with the 
property that $\pi^*\MF{M}\cdot r\cdot t \cdot b
\subseteq H^0(\CT,L)$.  We say that the decomposition satisfies $(\dag)$ if
$\MF{M}|_{C_i}\ne 0$ for all irreducible curves $C_i\subseteq C$.

We have the following lemma:
\begin{lem}\label{glb3}
If $H^0(L)|_{\CT_i}\ne 0$ for all $1\le i \le \mu$,
there is a decomposition of $L$ satisfying $(\dag)$.  Moreover, 
given a   decomposition of $L$ satisfying $(\dag)$ with $B=0$, there is an isomorphism
$\tt^*(\nu^*L)\cong \nu^*L$ such that
$$H^0(\NT,\nu^*L)^+=\pi^*H^0(N,M)\cdot r \cdot t,$$
where $H^0(\NT,\nu^*L)^+$ refers to the positive eigen-space of  the induced action of $\tau$ on $H^0(\NT,\nu^*L)$.
\QED
\end{lem}

\remskip
For irreducible curves, there are the following two results, which
closely reflect those in \cite{casa} for smooth curves.  The proofs are similar.  


\begin{cor}[\cite{casa}, Corollary 2.1.3]\label{glbirr}
If $\CT$ is irreducible  there exists a decomposition of
$L$ satisfying  $(\dag)$.  Moreover, given such a decomposition, let $n_1=\dim(\MF{M})$, 
let $n_2=h^0(L)-\dim(\MF{M})$, and let 
 $p_1,\ldots,p_k$ be general points of $\CT$. Set $D_k=\sum_{i=1}^k (p_i+\tt(p_i))$ and $\bar{p}_i=\pi(p_i)$.
\begin{itemize}

\item[\textnormal{(a)}] If $\dim(\MF{M})>h^0(L)/2$ and $k\le n_2$, $h^0(L(-D_k))=h^0(L)-2k$.

\item[\textnormal{(b)}] If $\dim(\MF{M})>h^0(L)/2$ and 
$n_2\le k\le n_1$,
then $h^0(L(-D_k))=h^0(L)-n_2-k$.  In this case,
$H^0(L(-D_k))=\pi^*\MF{M}(-\sum_{i=1}^k\bar{p}_i)\cdot r \cdot t\cdot b$.

\item[\textnormal{(c)}] If $\dim(\MF{M})>h^0(L)/2$, 
$n_2\le k\le n_1$, and $1\le k_1\le k$, 
$$h^0(L(-D_k-D_{k_1}))=h^0(L(-D_k-\sum_{i=1}^{k_1}p_i))=\max(h^0(L)-n_2-k-k_1,0).$$
In this case, $
H^0(L(-D_k-D_{k_1}))=\pi^*\MF{M}(-\sum_{i=1}^k\bar{p}_i-\sum_{i=1}^{k_1}\bar{p}_i)
\cdot r \cdot t\cdot b$.
\QED
\end{itemize}
\end{cor}

\begin{cor}[\cite{casa}, Lemma 2.1.4]\label{glbir}
Suppose  $\CT$ is irreducible, 
$h^0(L)=2n>0$, and for  general points
$p_1,\ldots,p_n$ of $\CT$,
$h^0(L(-D_n))>0$, where 
$D_n=\sum_{i=1}^n (p_i+\tt(p_i))$.
Then $L$ has a decomposition satisfying $(\dag)$ 
such that $\dim(\MF{M})>h^0(L)/2$.
\QED
\end{cor}

In general, for reducible curves it is hard to understand how points impose conditions on sections of line bundles. For curves of low genus we can say a little more.
The proof of the following proposition is similar to 
the arguments in \cite{cmf} on 
pages 27-30.  Recall the parity of 
a connected \'etale double cover of a plane quintic is 
that of $h^0(\OO_{\PR{2}}(1)|_C 
\otimes \eta)$, where $\eta$ is the two-torsion line
bundle associated to the double cover.

\begin{pro}\label{mainpro1} Suppose $p_a(C)\le 6$.
If $h^0(L)=4$, and 
there exists a decomposition of $L$ satisfying $(\dag)$, 
then $h^0(N,M)\le 3$.  If $\dim(\MF{M})=  3$, then 
\begin{itemize}
\item[\textnormal{(a)}]  if $p_a(C)=5$, then $C$ is hyperelliptic;

\item[\textnormal{(b)}] if $p_a(C)=6$,
then $C$ is either 
hyperelliptic, or obtained from a hyperelliptic curve by identifying two points,
or trigonal, or is a plane quintic, 
in which case $\pi$ is an odd double cover.\QED
\end{itemize}
\end{pro}

The following proposition complements Proposition \ref{mainpro1}.    The proof is similar to those in \cite{casa}.
\begin{pro}\label{mainpro2} Suppose $p_a(C)\le 6$,
$h^0(L)=4$, there exists a decomposition of
$L$ satisfying $(\dag)$, and for every such decomposition
$\dim(\MF{M})<3$.  There exist smooth points $p_1,p_2 \in \CT$ such
that $H^0(L(-\sum_{i=1}^2(p_i+\tt(p_i))))=0$. \QED
\end{pro}

\section{Prym theta divisors}

The results of Section 2 will allow us to determine the 
multiplicity of a point on the Prym theta divisor in certain cases.  

\subsection{General statements}

The following result was proven in \cite{cmf}:

\begin{pro}[\cite{cmf}, Lemma 2.1] \label{prst1}
Suppose $x\in P$ corresponds to a line bundle $L\in P^*$ such that
$h^0(L)=2n$, and there exist $n$ distinct points $p_1,\ldots, p_n$ of $\CT$ 
such that $H^0(L(-D_n))=0$, where 
$D_n=\sum_{i=1}^n (p_i+\tt(p_i))$.
Then $\MUL_x\Xi=h^0(L)/2$.
\QED
\end{pro}

The next proposition rephrases a result from \cite{cmf} in 
the notation 
used in Section 2.

\begin{pro}[\cite{cmf}, Theorem 2.3]\label{prst2}
If  $x\in P$ corresponds to a line bundle $L\in P^*$ 
which has a decomposition satisfying $(\dag)$,
$\MUL_x \Xi \ge \dim(\MF{M})$.
\end{pro}

\PRF.  The proof follows directly from \cite{cmf} Theorem 2.3, and Lemma
\ref{glb3}. \QED 


We also have the following:
\begin{teo}\label{t2}
Suppose 
$x$ is a point of $\Xi$, 
corresponding to a line bundle $L\in P^*$ such that
$h^0(L)=2$.  If there exists a general point $p$ on 
$\CT$ such that $h^0(L(-p-\tt(p)))=0$, then $x$ is a smooth point of $\Xi$.  
Otherwise, $\MUL_x\Xi=2$.
\end{teo}
The first statement is a special case of Proposition \ref{prst1}.  
The second statement is a consequence of the following:

\begin{pro}\label{p2}
In the notation above, suppose $\CT=\bigcup_{i=1}^d \CT_i$, where
the $\CT_i$ are connected.
If $h^0(L)=n\le d$, $\dim(H^0(L)|_{\CT_i})=1$ for all $i\in \{1,\ldots, n\}$, and $\dim(H^0(L)|_{\CT_j})=0$ for all $j\in \{n+1,\ldots, d\}$,
then $\MUL_x \Xi= h^0(L)$.
\end{pro}

\PRF. 
To begin, one can easily check that for any deformation $\LL$ which lies in $P$,
all sections of $L$ lift to first order. In other words, in the notation of Section \ref{secpre} and Lemma \ref{lemest}, taking $\LL$ general and setting $\LL_k$ to be the truncation of $\LL$ to order $k$, $\ell (H^0(\LL_1))=2h^0(L)$, 
and for some $k>0$,

\begin{equation}\label{eqnin}
h^0(L)=\frac{1}{2}\ell (H^0(\LL_1))\le \frac{1}{2}\ell (H^0(\LL_k))= \MUL_x\Xi.
\end{equation}

We now show that for
a particular deformation, there are no nontrivial sections lifting to 
second order.  
For $i\in \{1,\ldots,n\}$ let $p_i\in \CT_i$ be a general point.  Let $S\subset \mathbb C$ be a  neighborhood of the origin, and let $q_i:S
\to U_i\subset \CT$ be an isomorphism centered at $p_i$.
There is a deformation $\LL$ of $L$  given by 
$$
L\otimes \mathcal{O}_{\CT} \left( p_1-\tt(p_1)-q_1+\tt(q_1)\right) \otimes \ldots
\otimes \mathcal{O}_{\CT} \left( p_n-\tt(p_n)-q_n+\tt(q_n)\right)
$$
corresponding  to a morphism $f:S\to P$ sending $0$ to $x$.
For more details on these deformations see \cite{cmf} and \cite{casa}.
 Let $t$ be a local coordinate on $S$ centered at $0$, and  in the notation of Section \ref{secpre}, define $\LL_k$ on $\CT_k=\CT \times S_k$ to be the truncation of $\LL$ to order $k$.  Locally, for an open set $U\times S_k\subset \CT_k$,  a section of $H^0(\LL_k)$ can be written as $\sum_{i=0}^k \si^{(i)}t^i$ for some functions $\si^{(i)}\in \OO_{\CT}(U)$.

Observe that a section $s\in H^0(L)$ lifts to second order only if
$s|_{\CT_1}\in H^0(L|_{\CT_1})$ lifts to second order.
Now let $\CT_v=\cup_{i\ne 1}\CT_i$, and suppose that $s\in H^0(L)$ is such that
$s|_{\CT_1}\not\equiv 0$, and $s|_{\CT_v}\equiv 0$. 
Let $\CT_1\cap \CT_v
=\{ x_1,\ldots ,x_m\}$, let
$L'=L|_{\CT_1}\otimes \OO_{\CT_1}(-\sum_{i=1}^m x_i)$,  and consider the  deformation $\LL'$ of $L'$ induced by $\LL$, which in the notation above is given by 
$$
L'\otimes \mathcal{O}_{\CT_1} \left( p_1-\tt(p_1)-q_1+\tt(q_1)\right).
$$
 
One can prove inductively that if
$\sum_{i=0}^k \si^{(i)}t^i \in H^0(\LL_k)$, then 
$\si^{(i)}(x_j)=0$ for $1\le i \le k$, $1\le j \le m$.
In other words, the section $s$  lifts to second order as a section of 
$\LL_2$ only if 
$s|_{\CT_1}\in H^0(L')$ lifts to second order as a section of $\LL'_2$. 
Since
$H^0(L')=H^0(L)|_{\CT_1}$, it follows that $h^0(L')=1$.
One can check that 
$h^0(\CT_1, L'(p_1+\tt(p_1)))=2$, and that 
$h^0(\CT_1, L'(p_1+2\tt(p_1)))=2$, and thus
the natural inclusion $H^0(L'(p_1+\tt(p_1)))\subseteq H^0(L'(p_1+2\tt(p_1)))$
is an equality.
It follows from \cite{casa} Lemma 1.3.4 that 
 only the trivial section of 
$L'$ lifts to second order.

We now complete the proof. Let $s\in H^0(L)$ be such
that $s|_{\CT_1}\ne 0$.  Then $s$ lifts to order two only if 
$s|_{\CT_1}$ lifts to order two as a section of $\LL'_2$.  We have seen that
only the trivial section lifts to second order as a section of $\LL'_2$.  
By symmetry, only the zero section of $L$ lifts to second order; in other words
$H^0(\LL_2)=H^0(\LL_1)$.  From Lemma \ref{lemest} and the inequality \eqref{eqnin}, it follows that 
$h^0(L)\le \MUL_x\Xi\le \frac{1}{2}\ell( H^0(\LL_1))=h^0(L)$.\QED

\remskip
\noindent\emph{Proof of Theorem \ref{t2}}.  If there exist smooth
points 
$p,\tt(p)$ such that  $H^0(L(-p-\tt(p)))=0$, then
$\MUL_x \Xi=1$, by Proposition \ref{prst1}.  So suppose that
there do not exist such points.  
In the case where  there is no 
component $\CT_j\subseteq \CT$ such that $\dim (H^0(L)|_{\CT_j})=2$, 
the theorem follows from Proposition \ref{p2}.  
In the case where 
there is a 
component $\CT_j\subseteq \CT$ such that $\dim(H^0(L)|_{\CT_j})=2$, it is straightforward to adapt the proof of 
\cite{casa} Lemma 4.1.1 and \cite{casa} Theorem 2 to show that $\MUL_x\Xi=2$.
\QED

\comment{
In the case where 
there is a 
component $\CT_j\subseteq \CT$ such that $H^0(L)|_{\CT_j}=2$, 
let $\CT_{n}$ be the union of components $\CT_i$
such that $H^0(L)|_{\CT_i}\ne 0$, and let $s\in H^0(L)$ be such
that $s|_{\CT_i}\ne 0$ for all irreducible components of $\CT_n$.  
Then $s$ defines a divisor 
$(s)_0$ on $\CT_n$, such that $\OO_{\CT_n}((s)_0)=
L|_{\CT_n}$.  Let us call this line bundle $L'$.
$L'$ satisfies $(\dag)$, and clearly $\MF{M}'\cdot r \cdot t' \cdot b'
\cap H^0(L)\ne \varnothing$, due to $s$.  
I claim that this implies
$\MF{M}'\cdot r \cdot t \cdot b
\cap H^0(L) = H^0(L)$, because otherwise there would be smooth points
$p,\tt(p)$ imposing independent conditions on $H^0(L)$.  
Indeed, 
consider a component
$\CT_j$ where 
$H^0(L)|_{\CT_j}=2$.  If there does not exist $v\in 
H^0(L)|_{\CT_j}-(\pi^*\MF{M}'\cdot r_{0;n}\cdot b_n\cap H^0(L))_{\CT_j}$, 
then $\dim(\pi^*\MF{M}'\cdot r_{0;n}\cdot b_n\cap H^0(L) )=2$, 
and we are done.  Otherwise,
one shows that $v/b'$ can not be $\tt$-invariant on $\CT_j$, as in 
 Lemma \ref{glb3}, and then one can conclude in the same way.

The proof of Proposition \ref{prst2} shows that, 
$\MUL_x\Xi \ge 2$.  On the other hand, it is
clear that the numerical conditions of \cite{casa} Lemma 4.1.1 hold in this case,
and so $\MUL_x \Xi \le 2$.
\QED}

\remskip
The following theorem 
extends a theorem of Smith and Varley's \cite{sv}, and a theorem of the author's \cite{casa}, to the case of an 
irreducible curve with nodes.  We continue to use the notation 
of Section 1.

\begin{teo}
 Suppose $\CT$ is irreducible.  In the notation above, if
$x\in \Xi$ corresponds to a line bundle $L\in P^*$, and 
$C_x\THT$ is the tangent cone to $\THT$ at $x$, then the following are
equivalent:

\begin{itemize}
\item[\textnormal{(a)}]$T_xP\subseteq C_x \THT$;
\item[\textnormal{(b)}]$\MUL_x \Xi >h^0(L)/2$;
\item[\textnormal{(c)}] $L$ has a decomposition satisfying $(\dag)$, with
$\dim(\MF{M})>h^0(L)/2$.
\end{itemize}
Moreover, if (c) holds, 
$\MUL_x\Xi=\dim(\MF{M})$.
\end{teo}
\PRF.  The proof of the equivalence follows directly from Corollary \ref{glbir},
and Proposition \ref{prst2}.
In the case where  (c) holds,  using
 Corollary \ref{glbirr} 
  it is straightforward to adapt the proof of 
\cite{casa} Lemma 4.1.1, and \cite{casa} Theorem 2
 to show $\MUL_x\Xi=\dim(\MF{M})$.
\QED

\begin{cor}\label{corprym}
Suppose  $(P,\Xi)\in \AA_{g-1}$ is 
an indecomposable Prym variety associated to 
an admissible double cover of an irreducible curve.
For
$x\in \Xi$,  
$\MUL_x\Xi\le g/2=(\dim(P)+1)/2$.
\QED
\end{cor}

\subsection{Proof of 
Theorem \ref{teoel}, and 
Theorem
\ref{theorem}}

In this section
 we will only be concerned with indecomposable ppavs, and so we will always assume that 
the double cover satisfies 
$(\ast)$, 
and 

\begin{itemize}
\item[$(\sharp)$] if there is a decomposition
$\CT=\CT_1\cup \CT_2$, with $\CT_1\cap \CT_2$ finite, then 
$\#(\CT_1\cap \CT_2)$ is even and $\ge 4$. 
\end{itemize}
Otherwise, by \cite{bshot}, 
Theorem 5.4, and Lemma 4.11, $(P,\Xi)$ is either reducible or a Jacobian.

\begin{teo}\label{teoB}
Suppose 
$(P,\Xi)$ is indecomposable, 
and $x\in \Xi$. 
If $p_a(C)=5$, 
then $\MUL_x\Xi \le 2$.
\end{teo}

\PRF.  Suppose $x$ corresponds to the line bundle
$L\in P^*$. 
By Lemmas 3.15 and 3.16 of \cite{cmf}, we may assume $h^0(L)<6$.  If $h^0(L)=2$, then
by Theorem \ref{t2} $\MUL_x\Xi \le 2$.
If $h^0(L)=4$,  then again, by  Lemmas 3.15 and 3.16 of \cite{cmf}, $H^0(L)$ is nonzero on
every irreducible component of $\CT$, and so 
$L$ has a decomposition satisfying $(\dag)$.
If there does not exist such a decomposition with 
$\dim(\MF{M})\ge 3$, it follows from Proposition \ref{mainpro2} that for general points $p_1$ and $p_2$ of $C$,
$H^0(L(-p_1-\tau(p_1)-p_2-\tau(p_2)))=0$;  Proposition \ref{prst1} implies that in this case
$\MUL_x\Xi =2$.
If on the other hand there is a decomposition of $L$ with 
$\dim(\MF{M})\ge 3$, then by virtue of Proposition \ref{mainpro1},
$C$ is hyperelliptic.  A result of Mumford's implies that 
$(P,\Xi)$ is the Jacobian of a curve of genus four.  
The Riemann singularity theorem implies that $\MUL_x \Xi \le 2$. \QED

\remskip
\begin{rem}\label{remA}
It has been shown by Beauville in \cite{bshot} that if 
$\MUL_x\Xi=2$, then 
$(P,\Xi)$ is the Jacobian of a curve, or has a vanishing theta null.
\end{rem}

\begin{teo}\label{teoA}
Suppose 
$(P,\Xi)$ is indecomposable, 
and $x\in \Xi$. 
If $p_a(C)=6$, then $$\MUL_x \Xi \le 3.$$
Moreover, if equality holds
then $C$ is 
either hyperelliptic, or 
obtained from a hyperelliptic curve by identifying two points,
or trigonal, or a plane quintic, in which case $\pi:\CT\to
C$ is the double cover associated to an odd theta characteristic.
It follows that $(P,\Xi)$ is either the Jacobian of a 
hyperelliptic curve, or the
intermediate Jacobian of a smooth cubic threefold in $\PR{4}$.
\end{teo}

\PRF.  Suppose $x$ corresponds to the line bundle
$L\in P^*$, and $\MUL_x \Xi\ge 3$.
Again, by Lemmas 3.15 and 3.16 of \cite{cmf}, we may assume $h^0(L)\le 6$, and 
equality holds only if $C$ is hyperelliptic, and hence
$(P,\Xi)$ is the Jacobian of a hyperelliptic curve.
If $h^0(L)=2$, then
by Theorem \ref{t2} $\MUL_x\Xi \le 2$, which is a contradiction.
If $h^0(L)=4$,   
then $H^0(L)$ is nonzero on
every irreducible component of $\CT$, and so 
$L$ has a decomposition satisfying $(\dag)$.
We claim that  there exists such a decomposition of $L$ with 
$\dim(\MF{M})\ge 3$.  Indeed, if not,
 it would follow from Proposition \ref{mainpro2} that 
for general points $p_1$ and $p_2$ of $C$, 
$H^0(L(-p_1-\tau(p_1)-p_2-\tau(p_2)))=0$;  Proposition \ref{prst1} would then imply that $\MUL_x\Xi =2$, which is a contradiction.
Thus there is a decomposition of $L$ with 
$\dim(\MF{M})\ge 3$, and the claims on $C$
follow directly from
Proposition \ref{mainpro1}.
Results of Mumford \cite{mum}, Shokurov
\cite{shok}, Beauville \cite{bshot} and \cite{bdet}, [cf. 
\cite{cmf}, Theorem 4.1], then imply that
$(P,\Xi)$ is either the Jacobian of a curve, or the
intermediate Jacobian of a cubic threefold.  It follows that
$\MUL_x \Xi \le 3$, and equality holds only if
$(P,\Xi)$ is the Jacobian of a hyperelliptic curve, or the
intermediate Jacobian of a cubic threefold.
\QED

We have the following consequence.

\begin{cor}\label{teo1} An indecomposable ppav of  dimension five  whose theta divisor has a triple point is either the intermediate Jacobian of a smooth cubic threefold in $\PR{4}$, or the Jacobian of a hyperelliptic curve.  \end{cor} 

In \cite{cmf} the same statement was proven for ppavs of dimension five with a unique triple point.  

\remskip
\noindent\emph{Proof of Corollary \ref{teo1} and Theorem \ref{theorem}}.
The statements 
are now  a direct consequence of Theorem \ref{teoA},
Theorem \ref{teoB}, 
Corollary \ref{corprym}.\QED

\remskip
Theorem \ref{teoel} will follow easily from these results, together with the lemma below, which is essentially proven in
Ein and Lazarsfeld \cite{einlaz}.

\begin{lem}\label{lemel} Let $N$ and $k$ be a positive integers.
The following statements are equivalent:
\begin{enumerate}
\item[\textnormal{(a)}]  For a ppav $(A,\Theta)$ of dimension $d\le N$, 
$\dim(\textnormal{Sing}_k\Theta)> d-2k+1$
 implies that $(A,\Theta)$ is reducible.
\item[\textnormal{(b)}] For a ppav $(A,\Theta)$ of dimension $d\le N$,
and all nonnegative integers $j$  such that $k-j\ge1$, if $\dim(\textnormal{Sing}_k\Theta)= d-k-j$,
then
$(A,\Theta)=\prod_{i=1}^{k-j}(A_i,\Theta_i)$. \QED
 \end{enumerate}
\end{lem}

\remskip
\noindent\emph{Proof of Theorem \ref{teoel}}.
By virtue of the lemma, we need only show that for 
$d\le 5$, if $(A,\Theta)\in \mathcal{A}_5$, 
and $\Theta$ is irreducible,
then
$\dim(\textnormal{Sing}_k\Theta)\le d-2k+1$.
By Theorems \ref{teoB}, and \ref{teoA}, the pertinent cases are $d=5$, $k\le 3$, and $d=4$, $k\le2$.
If $d=5$ and $k=3$, then the result follows from Theorem \ref{teoA}; if 
$d=5$ and $k=2$, or $d=4$ and $k=3$, then the result follows from 
Ein and Lazarsfeld \cite{einlaz}, Corollary 2 .  
\QED

\remskip
\begin{rem}
For the Jacobian of a curve, Martens' theorem implies that 
$\dim(\singd_k\Theta)=g-2k+1$ only if the curve is hyperelliptic.  
It is a result of Beauville \cite{bshot} that if
$(A,\Theta)$ is an indecomposable Prym variety, and 
$\dim(\singd_2\Theta)\ge d-4+1$, then 
$(A,\Theta)$ is a hyperelliptic Jacobian.  
Thus at least in dimension five or less, any indecomposable ppav whose theta divisor has double points in codimension three is a hyperelliptic Jacobian.  In regards to these results, it was asked in \cite{casa} to what extent 
$k$-fold points in codimension $2k-1$ on an indecomposable ppav characterize hyperelliptic Jacobians.  Theorem \ref{teoB} implies that in dimension less than or equal to five, 
$k$-fold points in codimension $2k-1$ imply that $(A,\Theta)$ is either a hyperelliptic Jacobian, or the intermediate Jacobian of a smooth cubic threefold.
\end{rem}

\subsection{The triple point locus}\label{geomsec}

Recall  the definition of the Andreotti-Mayer 
loci $N_k\subseteq \AA_g$:

$$
N_k=\{(A,\Theta)\in \AA_g:\dim(\sing \Theta)\ge k\}.
$$
These are closed subschemes, and similarly, it is easy to see that the subloci 
$$
S_k=\{(A,\Theta)\in \AA_g: \singd_k \Theta \ne \varnothing\}
$$
are closed as well.

It is a result
of Andreotti and Mayer \cite{am} that 
$\bar{J}_g$ is an irreducible component of 
$N_{g-4}$.  
In \cite{bshot}, Beauville showed that
$N_0\subseteq \AA_4$ consists of two irreducible components of dimension 
nine; to be exact, $N_0= \bar{J_4}\cup \theta_{\textnormal{null}}$, 
where $\theta_{\textnormal{null}}$ is the locus of ppavs with a vanishing 
theta null.  We will prove a similar statement for $S_3\subset \mathcal A_5$.  

Let 
$I\subseteq \AA_5$
 denote the locus of intermediate Jacobians of smooth cubic threefolds
in $\PRRRR$,
and let $J^h_{g_1,\ldots,g_n} \subseteq \AA_{\sum_i g_i} $ be the
locus of ppavs which are the product of $n$
Jacobians of hyperelliptic curves of genera  $g_1,\ldots,g_n$.
For a sublocus $V\subseteq \AA_g$, let
$\AA_{g'}V\subseteq \AA_{g'+g}$ be the locus of ppavs which are the product of
a ppav of dimension $g'$ with  
a ppav in $V$.

\begin{teo}\label{teos3}  
The locus $S_3\subset \mathcal A_5$ has three irreducible components, $\bar{I}$, $\AA_1\bar{J}_4$ and $\AA_1\theta_{\textnormal{null}}$, each of which is ten-dimensional.
\end{teo}

\PRF.  Let $(A,\Theta)\in S_3$, and suppose first that 
$(A,\Theta)=(A_1,\Theta_1)\times (A_2,\Theta_2)$. 
We may as well assume that $\dim(A_1)<\dim(A_2)$. 
If $\dim(A_1)=1$, then $\Theta_2$ must have 
a double point, and so by Beauville's result on $N_0$, 
$(A,\Theta)\in
\AA_1\bar{J}_4\cup \AA_1 \theta_{\textnormal{null}}$. If $\dim(A_1)=2$, 
then we may assume that $(A_1,\Theta_1)$ and $(A_2,\Theta_2)$ are
indecomposable, since otherwise
we would be in the previous case.  It follows that $\Theta_2$  
has a double point, and so $(A_2,\Theta_2)$ is 
the Jacobian of a hyperelliptic curve.  
Thus, $(A,\Theta)\in J^{h}_{2,3}$.  A result
of Collino's \cite{collino} (cf. the remark below) states that
$J^h_5\cup J^h_{1,4}\cup J^h_{2,3}\subseteq \partial I$;
thus
$(A,\Theta)\in \bar{I}$.
Finally, if $(A,\Theta)$ is indecomposable, then it follows from Corollary
\ref{teo1}
that either $(A,\Theta)$ is the intermediate Jacobian of a smooth cubic 
threefold, or it is the Jacobian of a hyperelliptic curve.  Again,
Collino's result  implies that $(A,\Theta)\in \bar{I}$.  
The fact that $\AA_1\theta_{\textnormal{null}}$
is irreducible and of dimension ten follows from the fact (cf. Beauville \cite{bshot})
that $\theta_{\textnormal{null}}$ is irreducible of dimension nine.
\QED

\remskip
\begin{rem}\label{remsam}
Samuel Grushevsky has pointed out the following observation: for a ppav $(A,\Theta)\in \AA_g$ and a point $x\in A$ such that
$x+x=0$ and $\MUL_x\Theta=3$, the gradient of the theta function at $x$ vanishes, giving $g$ conditions in $\AA_g$ cutting out $S_3$.  
The fact that in $\AA_5$, $S_3$ is equidimensional and has codimension five 
is equivalent to the independence of the conditions arising from the vanishing gradient. 
Note also that in $\AA_4$, $S_3=J^h_{1,3}$, 
which has codimension four, 
and in $\AA_3$, $S_3=J^h_{1,1,1}$ which has codimension three.
Hence the conditions arising from the vanishing gradient are independent in these cases as well.  It would be interesting to know if this were the case for all $S_3\subset\AA_g$ with $g\ge 3$.
\end{rem}

\begin{cor}[Theorem \ref{teo0}]
In the notation above, $I=S_3-(N_1 \cap S_3)$.
\end{cor}
\PRF. Theorems \ref{teoA} and \ref{teos3} imply that 
the ppavs in $S_3$ which are not in $I$ are either reducible, or are
Jacobians of  hyperelliptic curves, and therefore the dimension of the
singular locus of their theta divisor must be at least one (in fact at least two).  \QED

\remskip

\begin{rem}
The fact that 
  $J^h_5\cup J^h_{1,4}\cup J^h_{2,3}
\subseteq \partial I$ was proven by Collino in 
\cite{collino} by considering the secant variety to the 
rational normal quartic.  It is interesting to point out
that this result can be obtained as well 
by combining some arguments used in Beauville \cite{b2}, Section 2.9,
with a result of Griffin's \cite{griff} on degenerations of plane
quintics, and a result of Mumford's \cite{mum} on Prym varieties of 
hyperelliptic curves.  
\end{rem}

\remskip
\begin{rem}\label{rembound}
On the other hand, it is still not clear exactly which ppavs lie in $\partial I$.  The results of this paper tell us that the only indecomposable ppavs in $\partial I$ are the hyperelliptic Jacobians.  In general, the description of $S_3$ above certainly narrows the list of possibilities, and one may be able to use similar techniques to give a complete description.

One can also consider the boundary of $I$ in the
Satake compactification, and 
certain special cases are known; e.g.
in \cite{cg}, Clemens and Griffiths show that
 the intermediate Jacobian of a cubic hypersurface in $\mathbb P^4$ with a unique ordinary double point  is 
the extension by a $\CX^*$ of the Jacobian of a curve of genus $4$.  
Other treatments and examples are considered by 
Collino and Murre \cite{collinomurre}, and 
Gwena \cite{gwena}.
\end{rem}

\remskip
\begin{rem}
In this paper, we study the locus of intermediate Jacobians of  cubic threefolds by classifying those ppavs of dimension five whose theta divisor has a triple point.
It should be pointed out that there are other properties of the theta divisor which have been conjectured to classify these intermediate Jacobians as well.   

For example, for a cubic threefold $X$, the Fano surface $F$ of lines on $X$ embeds in $JX$ as a nondegenerate surface, and has class $[F]=[\Theta]^3/3!\in H^6(JX,\mathbb{Z})$.   In \cite{deb}, Debarre showed that $\bar{I}$ is an irreducible component of the locus of ppavs of dimension five for which there exists a surface $S\subseteq X$ such that $[S]=[\Theta]^3/3!$.  
For more discussion in this vein, we refer the reader to Debarre \cite{deb}.

In another direction, Pareschi and Popa \cite{pp} have considered a notion of regularity for subvarieties of  abelian varieties which may also give a criterion for a ppav of dimension five to be the intermediate Jacobian of a cubic threefold (cf. Pareschi and Popa \cite{pp} Conjecture 2.2).

\end{rem}

\bibliography{ppavbib}

Harvard University

Department of Mathematics

One Oxford Street

Cambridge, MA 02138
\remskip

\texttt{casa@math.harvard.edu}

\remskip

Department of Mathematics

SUNY Stony Brook

Stony Brook, NY 11794-3651

\remskip
\texttt{casa@math.sunysb.edu}

\end{document}